\documentclass[12pt]{amsart}
\usepackage[francais]{babel}
\usepackage{graphics}
\usepackage{graphicx}
\usepackage{pdfpages}

\title{Les groupes de Mathieu sont-ils aussi sporadiques
qu'il y para\^it ?} 

\author{Labib Haddad}
\address{120 rue de Charonne, 75011 Paris, France}
\email{labib.haddad@wanadoo.fr}

\usepackage{amssymb}
\usepackage{amsmath}
\usepackage{endnotes}

\newcommand{\su}{\subsection*}
\newcommand{\head}{\section*}
\newcommand{\noi}{\noindent}
\newcommand{\se}{\noi{\bf En effet}}
\newcommand{\cqfd}{\hfill{\bf cqfd}}

\newcommand{\Ž}{\'e}

\newcommand{\ˆ}{\`a}
\newcommand{\}{\`u}

\newcommand{\ž}{\^u}

\newcommand{\Z}{\mathbb Z}

\newcommand{\cal}{\mathcal}

\newcommand{\leqs}{\leqslant}

\newcommand{\geqs}{\geqslant}

\newcommand{\guil}{\guillemotleft}  
\newcommand {\guir}{\guillemotright}

\newcommand {\et}{\ \text{et}\ }

\newcommand {\donc}{\ \text{donc}\ }

\newcommand {\pour}{\ \text{pour}\ }

\newcommand {\pourtout}{\ \text{pour tout}\ }

\newcommand{\stm}{\smallsetminus}

\newcommand{\vide}{\emptyset}

\newcommand{\inc}{\subset}

\newcommand{\bc}{\begin{cases}}
\newcommand{\ec}{\end{cases}}

\newcommand{\h}{h^{-1}}

\begin{document}
\maketitle

\thispagestyle{empty}

\markboth{LH}{Groupes de Mathieu sporadiques ?}

\

\

\

\hfill{\it \`A Yves Sureau, guide avis\Ž, de haute
}

\hfill{\it montagne, au pays des hypergroupes.}

\

\

\su{Abstract}{\sl This is a plea for the reopening of the building site for the
classification of finite simple groups in order to include
the finite simple hypergroups.

\

\noi Hypergroups were first introduced by Fr\Žd\Žric Marty, in 1934, at a Congress in Stockholm,  not to be confused with a later and quite different notion to which the same name was given, inopportunely.

\

\noi I am well aware that, probably,  quite a few mathematicians must have already felt uncomfortable about the presence of the so-called {\it sporadic simple groups} in the large {\rm tableau} of the classification of finite simple  groups, and might have wrote about it, though I do not have any reference to mention.

\

\noi In what follows, I will try to explain, step by step, what a hypergroup is, and, then, suggest a  notion of simplicity for hypergroups, in a simple and natural way,  to match the notion in the case of groups, hoping it will be fruitful.

\

\noi Examples and constructions are included.}

\newpage

\

\head{Introduction} 

\

Ce qui suit est un plaidoyer pour la r\Žouverture du chantier de
classification des groupes simples finis afin d'y inclure
les hypergroupes simples finis.

\

\noi La notion d'hypergroupe a \Žt\Ž introduite par Fr\Žd\Žric Marty, en 1934, lors d'un congr\s \ˆ Stockholm, \ˆ ne pas confondre avec une notion plus tardive et tout \ˆ fait diff\Žrente \ˆ laquelle on a, malencontreusement, donn\Ž le m\me nom.

\

\noi Je suis bien conscient qu'un bon nombre de math\Žmaticiens ont d\Žj\ˆ d\ž se sentir mal \ˆ l'aise au sujet de la pr\Žsence de ce que l'on nomme {\it les groupes simples sporadiques} dans le grand tableau de la classification, et certains ont  m\me pu le coucher par \Žcrit, bien que je ne puisse mentionner aucune r\Žf\Žrence.

\

\noi Dans ce qui suit, on va commencer par rappeller,  pas \ˆ pas,  ce qu'est un hypergroupe puis comment y introduire la simplicit\Ž de mani\re simple (et utile) et qui co\•ncide avec la notion classique, dans le cas des groupes. 

\

\noi On \Žlabore en donnant une mani\re tr\s g\Žn\Žrale pour la construction d'hypergroupes et l'on donne des conditions pour qu'ils soient simples. On fournit de nombreux exemples et caract\Žrisations

\

\noi On pose enfin la question : \guil \`A peine clos, le chantier de la classification des
{\bf groupes} simples finis devra-t-il, \Žventuellement,  rouvrir ses portes
afin d'entreprendre la classification des {\bf hypergroupes} simples finis?\guir

\

\head{Quelques d\Žfinitions} 

\

Un hypergroupe est un ensemble muni d'une op\Žration
binaire  associative {\bf multivalente} dont les propri\Žt\Žs g\Žn\Žra\-lisent celles des groupes. 

Plus pr\Žcis\Žment, une {\bf
op\Žration binaire multivalente} sur un ensemble donn\Ž
$H$ est une application $(x,y) \mapsto x.y$ qui, \ˆ
chaque couple d'\Žl\Žments $(x,y) \in H \times H$, fait
correspondre une {\bf partie} $x.y \inc H$. 

Pour 
$X \inc H$ et $Y \inc H$, on d\Žsigne alors par
$X.Y$ la r\Žunion de toutes les parties $x.y$ o\ $x$
parcourt
$X$ et
$y$ parcourt
$Y$.

Le maniement des op\Žrations multivalentes n\Žcesssite un
tout petit effort d'attention de plus que l'usage des op\Žrations classiques, univalentes, n'en demande, mais l'habitude finit par
estomper la difficult\Ž.

Lorsque le risque de confusion est minime, on
\Žcrit, simplement,  $x.Y$ au lieu de $\{x\}.Y$ et, de m\me,
$X.y$ au lieu de $X.\{y\}$. Enfin, si une des parties
$x.y$ est un singleton, au lieu de $x.y = \{z\}$, on \Žcrit
plus simplement $x.y = z$ et on retrouve l'\Žcriture classique des
op\Žrations binaires {\it univalentes}.

\su{D\Žfinition} Un {\bf hypergroupe} est un ensemble non vide, $H$, muni
d'une op\Žration binaire multivalente qui v\Žrifie les deux
identit\Žs suivantes, pour tous
$x,y,z$ dans $H$,

\

$x.(y.z) = (x.y).z$ : associativit\Ž.

$x.H = H = H.x$ : reproductivit\Ž.

\

 Ces deux propri\Žt\Žs impliquent que $(x.y).H = x.(y.H) =
x.H = H$, de sorte que le produit $x.y$ n'est jamais vide.

Bien entendu, tout groupe est un hypergroupe dont l'op\Žration est univalente. R\Žcipro\-quement, tout hypergroupe
dont l'op\Žration est univalente est un groupe
puisqu'alors, pour chaque couple $(a,b)$ donn\Ž, les \Žquations 
$a.x = b \et  y.a = b$ poss\dent des solutions en $x$ et en $y$.

\su{Un exemple \Žl\Žmentaire} Soient $G$ un groupe, $H$ l'un quelconque de ses
sous-groupes  et soit $G/H = \{xH : x \in G\}$ l'ensemble des {\it
classes  \ˆ \ droite} de $G$ {\it modulo} $H$.

L'ensemble $xHyH$ est la r\Žunion des classes \ˆ 
droite $xhyH$ o\ \ $h$ parcourt $H$. En posant
$(xH).(yH) = \{xhyH : h \in H\}$, 
on d\Žfinit une op\Žration binaire multivalente sur l'ensemble $G/H$ qui en fait
un hypergroupe, comme on le v\Žrifie sans d\Žtour. 

Cette op\Žration est univalente si et seulement si $H$ est un
sous-groupe invariant de $G$ : dans ce cas l'hypergroupe
$G/H$ est le groupe quotient classique de $G$ par $H$.
Lorsque le sous-groupe $H$ n'est pas invariant dans
$G$, on obtient alors un exemple {\it \Žl\Žmentaire}
d'un hypergroupe qui n'est pas un groupe.

Traditionnellement, on appelle {\bf D-hypergroupe} tout
hypergroupe isomorphe \ˆ un hypergroupe de classes \ˆ 
droite de la forme $G/H$. 

\su{Morphismes}  Un {\bf morphisme} est alors, par d\Žfinition, une application $f : H \to K$
d'un hypergroupe $H$ dans un hypergroupe $K$ qui v\Žrifie
l'identit\Ž d'homomorphisme, $f(x.y) = f(x).f(y)$.

\

On appelle 
{\bf isomorphisme} entre $H$ et $K$ toute  application
{\bf bijective} $f : H \to K$ qui est un morphisme ainsi
que sa r\Žciproque
$f^{-1}$. Comme dans le cas des groupes, un morphisme
bijectif $f : H \to K$ est (d\Žj\ˆ) un isomorphisme  car,
dans ce cas, $f^{-1}$ est aussi un morphisme.  En effet, pour
$a = f(x)$ et
$b = f(y)$, on aura
$x.y = f^{-1}f(x.y) = f^{-1}(f(x).f(y))$, autrement dit, 
$f^{-1}(a).f^{-1}(b) = f^{-1}(a.b)$.

Bien entendu, lorsque $H$ et $K$ sont des groupes, on
retrouve ainsi les notions classiques d'homomorphismes et
d'isomorphismes de groupes.

\su{Nombres premiers et groupes simples } Comme pour la primalit\Ž, il y a deux mani\res \Žquivalentes de d\Žfinir la simplicit\Ž.

\

Un entier $n \geqs 0$ est {\bf premier} lorsqu'il a
exactement deux diviseurs : cela  revient \ˆ dire qu'il
est seulement multiple de deux entiers distincts.

\

Un groupe est {\bf simple} lorsqu'il a exactement
deux sous-groupes invariants : cela revient \ˆ dire qu'il
poss\de seulement deux images homomorphes distinctes, \ˆ 
isomorphisme pr\s.

\

Afin d'\Žtendre aux hypergroupes la notion de simplicit\Ž,
il faut trouver un substitut convenable
\ˆ  la notion {\it d'image homomorphe} \rm d'un groupe. Voici ce que nous proposons.

\su{Reflets et simplicit\Ž} On  dira que l'hypergroupe 
$K$ est un {\bf reflet}
de l'hypergroupe $H$ lorsqu'il existe une application
{\bf surjective}
$f : H \to K$ qui v\Žrifie les identit\Žs suivante :

\

$f^{-1}f(x.y) = f^{-1}(f(x).f(y)) = x.f^{-1} f(y)
=f^{-1}f(x).y$.

\

On dira alors que l'application
$f$ est un
{\bf r\Žflecteur}. 

\

Tout r\Žflecteur $f$ est n\Žcessairement un morphisme. En
effet, $f$ \Žtant surjective, on a $f(x.y) = ff^{-1}f(x.y)
= ff^{-1}(f(x).f(y)) = f(x).f(y)$.

\

On observera aussi, sans d\Žtour, qu'une application $f : H
\to K$ est un isomorphisme si et seulement si c'est un r\Žflecteur injectif.

\

Bien entendu, si $H$ est un groupe, l'hypergroupe $K$ est un
reflet du groupe $H$ si et seulement si c'est une image homomorphe
de $H$.

\

Chaque  hypergroupe $H$ poss\de toujours, naturellement,
les deux reflets suivants : lui-m\me et l'hypergroupe {\bf
trivial}
 \ˆ un seul
\Žl\Žment; ces deux reflets sont isomorphes si et seulement
si $H$ lui-m\me est trivial.

\

On dira qu'un hypergroupe est {\bf simple} lorsqu'il
poss\de (\ˆ isomorphisme pr\s) exactement deux
reflets. Cette mani\re de proc\Žder permet d'\Žviter les
{\it arcanes} de la th\Žorie des sous-hypergroupes. On peut esp\Žrer que cette nouvelle notion soit un bon choix.

\

Il est clair qu'un groupe donn\Ž est simple si et seulement
si c'est un \guil hypergroupe simple\guir. La classe des
hypergroupes simples finis contient ainsi celle des groupes
simples finis. Elle la contient, {\bf et la prolonge} :
voir, en effet,  la cons\Žquence (5) ci-dessous. 

\su{Invariance} On aura besoin de la d\Žfinition suivante. Etant donn\Žs
un groupe $G$ ainsi que deux sous-groupes $H$ et
$K$, on dira que le sous-groupe $K$ est {\bf invariant modulo
$H$} lorsque l'on a
\[KxK = HxK = KxH, \pourtout x\in G.\tag {$*$}\]

Cela g\Žn\Žralise la notion classique de sous-groupe
invariant car, bien entendu, un sous-groupe est invariant
(au sens classique) si et seulement s'il est invariant modulo
le sous-groupe trivial.

On observera, en passant aux inverses, que chacune des deux
conditions suivantes est \Žquivalente \ˆ l'autre :
$$KxK = HxK, \pourtout
\ x
\in G,$$
$$KxK = KxH, \pourtout
\ x
\in G.$$

Chacune d'elles est  ainsi \Žquivalente \ˆ la condition
($*$). En outre, lorsque $K$ est
invariant modulo $H$, on a n\Žcessairement $H \subset K$
puisque, pour $x \in K$, la condition $HxK = KxK$ entra\"ne
$HK =K$. La condition ($*$) peut donc s'\Žnoncer \Žgalement comme ceci :
\[H \subset K \et Kx \inc HxK, \pourtout
\ x
\in G.\tag{$*$}\]

\

\su{Th\Žor\me} \sl Soient $G$ un groupe, $H$ et $K$ deux
quelconques de ses sous-groupes. Si $K$ est invariant modulo
$H$, le D-hypergroupe $G/K$ est un reflet du D-hypergroupe
$G/H$. R\Žciproquement, tout reflet de
$G/H$ est isomorphe \ˆ  l'un de ces D-hypergroupes $G/K$
\rm o\ $K$ est invariant modulo $H$.

\

Afin de ne pas couper le fil de l'expos\Ž, nous renvoyons
la d\Žmonstration de ce th\Žor\me \ˆ  l'\sc Appendice
\rm ci-dessous.

\

\head{Cons\Žquences}

\

\su{1} Tout reflet d'un D-hypergroupe est donc un
D-hypergroupe.

\su{2} Le D-hypergroupe $G/H$ est simple si et seulement si
$H$ est distinct de $G$  et que les seuls sous-groupes $K$
de $G$ invariants
 modulo $H$ sont
$G$ et
$H$ lui-m\me.

\su{3} En particulier, le D-hypergroupe $G/H$ est simple d\s lors que
$H$ est un sous-groupe \bf maximal \rm dans
$G$.

\su{4} Ainsi, lorsque
$H$ est un sous-groupe maximal du groupe
$G$ et qu'il n'est {\bf pas invariant} dans $G$,
l'hypergroupe
$G/H$ est simple et n'est {\bf pas un groupe}.

\su{5}  \bf Il existe donc des hypergroupes simples finis qui
ne sont pas des groupes\rm.

\

\head{Classes \ˆ gauche} 

\

On d\Žfinit, semblablement, l'hypergroupe $H\backslash  G =
\{Hx : x \in G\}$ des classes
\ˆ
gauche modulo $H$. Les hypergroupes isomorphes \ˆ ces
hypergroupes de classes \ˆ
gauche poss\dent les m\mes propri\Žt\Žs, corr\Žlatives, que les D-hypergroupes. Par exemple,
$H\backslash G$ est un groupe si et seulement si
$H$ est invariant dans
$G$ et, s'il en est ainsi, les deux groupes $G/H$ et
$H\backslash G$ co\•ncident, bien entendu. De m\me, \Žtant donn\Žs deux sous-groupes $H \inc K \inc G$,
l'hypergroupe
$K\backslash G$ est un reflet de $H\backslash G$ si $K$ est
invariant modulo $H$. 

Le th\Žor\me et ses cons\Žquences
se transposent, corr\Žlativement. En particulier
$H\backslash G$ est simple si et seulement si $G/H$ l'est.

\

\head{Hypergroupe oppos\Ž}

\

Soit $H$ un hypergroupe  avec son op\Žration $(x,y) \mapsto
x.y$. Comme dans le cas des op\Žrations univalentes, on d\Žfinit l'op\Žration {\bf oppos\Že}
$(x,y) \mapsto x \circ y = y.x$. Muni de cette nouvelle
op\Žration, $H$ est de nouveau un hypergroupe que l'on d\Žsigne par $H^\circ$ et que l'on appelle {\bf l'hypergoupe
oppos\Ž}.

\ 

Que dire de l'hypergroupe $(G/H)^\circ$ oppos\Ž d'un
hypergroupe
$G/H$?

\su{Th\Žor\me} \sl Soient $G$ un groupe et $H$ un de ses
sous-groupes\rm.

(1) \sl Les deux hypergroupes $(G/H)^\circ$ et
$H\backslash G$ sont isomorphes\rm.

(2) \sl D'un autre c\™t\Ž , pour que les deux hypergroupes
$G/H$ et $H\backslash G$ soient isomorphes, il faut et il
suffit que $H$ soit un sous-groupe invariant dans $G$\rm.

\

Autrement dit, lorsque le sous-groupe $H$ n'est pas
invariant, les deux hypergroupes $G/H$ et
$H\backslash G$ ne sont \bf pas isomorphes\rm. Lorsque le
sous-groupe $H$ est invariant, $G/H$ et $H\backslash G$ sont
un seul et m\me {\bf groupe} quotient.

\

\su{D\Žmonstration} (1) Plus pr\Žcis\Žment, l'application
$Hx \mapsto (Hx)^{-1}$ qui,
\ˆ
\ chaque classe \ˆ droite
$Hx
\in H\backslash G$, fait correspondre la classe \ˆ gauche
$(Hx)^{-1} = x^{-1}H
\in G/H$ est un isomorphisme de l'hypergroupe $H\backslash
G$ sur l'hypergroupe oppos\Ž \ $(G/H)^\circ$\rm, comme une
simple v\Žrification le montre.

(2) Si $H$ est invariant, alors $G/H$ et $H\backslash G$
sont, tous deux, le m\me groupe quotient classique de
$G$ par $H$.

R\Žciproquement, soit $f : G/H \to H\backslash G$ un
isomorphisme entre ces deux hypergroupes. Dans chacun d'eux,
$H$ est la seule classe qui v\Žrifie l'identit\Ž $H.H =
H$. On a donc $f(H) = H$. De plus, quelque soit $x \in G$, on
a
$(xH).H = xH$ dans l'hypergroupe $G/H$. Pour chaque $y
\in G$, il existe $x \in G$ tel que
$f(xH) = Hy$. On aura alors

\

$Hy = f(xH) = f((xH).H) = f(xH).f(H) = (Hy).H$.  

\

Ainsi, dans
$H\backslash G$, le produit  $(Hy).H = Hy$ est un singleton
: cela revient \ˆ  dire que $HyH = Hy$, de sorte que
$Hy \subset yH$ \bf pour tout \rm $y \in G$, donc
$H$ est invariant.\qed 

\

\head{Ils vont par paires} 

\

Soit $K$ un sous-groupe d'un groupe donn\Ž 
$G$. Si l'hypergroupe $G/K$ est simple, on sait que
l'hypergroupe $K\backslash G$ est \Žgalement simple.
Lorsque, de plus, ce ne sont pas des groupes, on sait qu'ils
ne sont pas isomorphes. Ainsi, {\bf comme les racines imaginaires
conjugu\Žes d'un polyn\™me r\Žel, ces hyergroupes simples
vont par paires !}

\

En particulier, en prenant tous les couples $(G,K)$ o\ \
$K$ est un sous-groupe maximal, non invariant, d'indice 
\bf fini\rm, dans un groupe
$G$, on obtient une famille de paires d'hypergroupes simples
finis $(G/K,K\backslash G)$, deux \ˆ \ deux oppos\Žs et non
isomorphes : cette famille
{\bf prolonge} ainsi la famille des {\bf groupes cycliques simples
finis}.

\

\head{Les plus simples des exemples} 

\

On prend pour $G$ le groupe de toutes les permutations d'un
ensemble $E$ non vide, fini ou infini, ayant $\alpha$ pour
cardinal. On distingue un point particulier
$p
\in E$ et on prend 
$H =
\{s
\in G : s(p) = p\}$, le sous-groupe des permutations qui
laissent le point
$p$ fixe. C'est un sous-groupe maximal de $G$.  Pour
$\alpha \geqs 3$, ce sous-groupe $H$ n'est pas invariant, de
sorte que
$K = G/H$ est alors un D-hypergroupe simple qui n'est \bf
pas un groupe \rm et dont le cardinal est celui de
$E$. En d\Žsignant par
$e$ la classe
$H \in K$, on  v\Žrifie ais\Žment que sa table de
multiplication se r\Žsume alors ainsi :
$$x.e = x, \pourtout x \in K,$$
$$x.y = K \setminus \{x\}, \pour \  y \neq e.$$

On observera encore ceci : pour tout
$x
\neq e$, on a
$$x^2 = x.x = K \setminus \{x\} \et \  x^3 =x.x.x = K.$$ 

Ainsi, contrairement au cas des groupes, {\it l'ordre} de chacun de ces hypergroupes
{\it monog\nes} est toujours \Žgal \ˆ {\bf trois}
et ne co\•ncide avec son cardinal que pour
$\alpha = 3$.

\

En particulier, en prenant pour $\alpha$, successivement, chacun des entiers
$n
\geqs 3$, on obtient une suite de D-hypergroupes simples
finis ayant $n$ \Žl\Žments et qui ne sont {\bf pas des
groupes}. 

\

\head{Un pas de plus}

\

Le proc\Žd\Ž de construction des hypergroupes de classes, \ˆ droite ou \ˆ gauche, 
de la forme
$G/H$ ou
$H\backslash G$, peut \tre  grandement \Žtendu en une {\bf
large g\Žn\Žralisation} que voici.

\

On se donne un ensemble $T$, une partie $C \inc T
\times T$ ainsi qu'une op\Žration binaire univalente,
$\rm{Op} : C \to T$, \bf partiellement \rm d\Žfinie,  pour
les seuls couples
$(x,y) \in C$ que l'on appellera les couples
{\bf composables}.  Le plus souvent, on \Žcrira simplement
$xy$ au lieu de
$\rm{Op}(x,y)$ et on dira que
$(T,C)$ est une
{\bf trame}.

\

Soit $R$  une relation {\bf d'\Žquivalence} sur l'ensemble $T$ et, pour chaque
\Žl\Žment $x
\in T$, soit $\bar x$ sa classe modulo $R$. On d\Žfinit
alors une op\Žration multivalente (naturelle) sur
l'ensemble des classes $T/R$, de la mani\re suivante :
$$\bar z \in \bar x .\bar y \iff (\exists (u,v) \in C)  (uv \in \bar z).$$

On dira que $(T,C,R)$ est une {\bf pr\Žsentation} de la
{\bf structure quotient}
$T/R$ et, lorsque
$T/R$ est un hypergroupe, on dira que la relation
$R$ est
{\bf ad\Žquate}.

\

On peut donner, facilement, des conditions n\Žcessaires et
suffisantes pour que la relation d'\Žquivalence $R$ soit
ad\Žquate, sous la forme d'un ensemble d'\Žnonc\Žs {\it du
premier ordre}. Il n'est pas besoin de rentrer dans le d\Žtail, ici. Voir, pour cela, l'\sc Appendice\rm \ ci-dessous. Pour le
moment, il suffit de savoir que ces conditions explicites
existent.

\

\head{Une heureuse circonstance}

\

Bien entendu, chaque hypergroupe de classes, de la forme
$G/H$ ou
$H\backslash G$, poss\de une pr\Žsentation {\it naturelle} dont la trame est $(T,C)$ o\ $T = G$ est le groupe
lui-m\me, $C = G \times G$, (i.e., tous les couples sont
composables, l'op\Žration \Žtant celle du groupe). La
relation d'\Žquivalence $R$ est d\Žfinie par
$$x R y \iff x^{-1}y \in H,$$ respectivement,
$$x R y \iff yx^{-1} \in H.$$

Plus g\Žn\Žralement donc, chaque hypergroupe isomorphe \ˆ
 l'un quelconque de ces hypergroupes de classes poss\de
\Žgalement une pr\Ž\-sentation. Cette particularit\Ž n'est,
cependant, pas due \ˆ
une quelconque singularit\Ž  de ces hypergroupes. Elle
est le lot de tous les hypergroupes, quels qu'ils soient.
Voici comment.

\

\centerline{\bf Tout hypergroupe poss\de une pr\Žsentation}

\

\noi \bf En effet\rm, soit $H$ un hypergroupe quelconque et
$T = H
\times H \times H \times H$. \`A chaque triplet $t =
(x,y,z)
\in H \times H \times H$ pour lequel $z \in x.y$, on associe
le couple \sl composable \rm $((x,t),(y,t)) \in T
\times T$ et on pose
$\rm{Op} ((x,t),(y,t)) = (z,t)$ : cela d\Žfinit une op\Žration
$\rm{Op} : C \to T$  o\ 
$C$ est l'ensemble de tous les couples composables. On
obtient ainsi la trame $(T,C)$. On prend pour $R$ la relation
d'\Žquivalence sur $T$ dont les classes sont les parties 
$\{x\}
\times H \times H \times H$. On v\Žrifie, sans d\Žtour,
que $T/R$ est un hypergroupe isomorphe \ˆ  l'hypergroupe
donn\Ž $H$.\cqfd

\

\head{Une nouvelle famille d'exemples}

\

Voici une premi\re illustration du proc\Žd\Ž de pr\Žsentation qui fournit des hypergroupes qui ne sont pas des D-hypergroupes.

\

On prend un ensemble $K$ r\Žunion d'une famille
$(A_i)_{i \in I}$ d'ensembles non vides, deux
\ˆ
deux disjoints. On suppose que $0 \in I$ et on distingue
un point particulier
$e
\in A_0$. On introduit l'ensemble $T$ form\Ž  de toutes les
applications
{\bf injectives} $f : X \to K$ o\  $e \in X \subset K$
et v\Žrifiant la condition suivante : il existe (au moins)
une permutation $s : I \to I$ telle que l'on ait, pour tout
$i \in I$,
$$f(X \cap A_i) \subset A_{s(i)}.$$

On dira que le couple $(f,g)$ form\Ž de deux \Žl\Žments $f : X \to
K$ et
$g : Y
\to K$ de $T$ est composable lorsque l'on a $g(e) \in X$ et
on prend pour compos\Ž $fg$ l'application $y \mapsto
f(g(y))$ d\Žfinie sur la partie $Y
\cap g^{-1}(X)$ : cette application, on le voit
sans d\Žtour, appartient bien \ˆ \ $T$. 

\

On obtient, ainsi,  une trame
$(T,C)$. L'application $f
\mapsto f(e)$ de
$T$ dans $K$ est surjective et d\Žfinit une relation d'\Žquivalence
$R$ sur $T$ dont les classes d'\Žquivalences sont les
parties $\bar x = \{f \in T : f(e) = x\}$. 

\

L'application
$x \mapsto \bar x$ permet d'identifier l'ensemble $T/R$ des
classes modulo $R$ \ˆ l'ensemble $K$ lui-m\me. L'op\Žration multivalente correspondant \ˆ  la pr\Žsentation
$(T,C,R)$ de $T/R = K$ est donn\Že par
$$z \in x.y \iff (\exists(f,g) \in C)(f(e) = x \ , \ g(e) =
y \ , \ f(y) = z).$$

On en \Žtablit facilement la table de multiplication sous la
forme suivante : quels que soient $a_i \in A_i$, pour chaque
indice $i$,

\

$x.e = x, \pourtout  x \in K$,

$a_i.y = A_i \setminus \{a_i\}, \pourtout  \ y
\in A_0 \setminus \{e\}$,

$a_i.a_j = K \setminus A_i, \pourtout \text{indice} \ 
j \neq 0$.

\

Dans le cas particulier o\  $I = \{0\}$, $K = A_0$, on
reconna\"t l'exemple du D-hypergroupe simple d\Žj\ˆ 
rencontr\Ž plus haut. 

\

Plus g\Žn\Žralement, si toutes les parties $A_i$ ont m\me cardinal, 
$K$ est encore un D-hypergroupe.

\

\se, on le voit facilement en introduisant
une nouvelle trame $G$, une sous-trame de la trame de d\Žpart $T$.  Du fait que les parties
$A_i$ sont
\Žquipotentes, chaque \Žl\Žment 
$f$ de
$T$ se prolonge en un \Žl\Žment $g$ de $T$ qui est une
application bijective
$g : K
\to K$, autrement dit, une permutation de l'ensemble $K$. Ce
prolongement n'est cependant pas n\Žcessairement unique. L'ensemble
$G
\inc T$ form\Ž  de toutes ces 
{\bf permutations} de $K$ est un groupe. On d\Žsigne par
$S$ la restriction de l'\Žquivalence $R$ au sous-ensemble
$G$ et par
$H$ le sous-groupe de $G$ form\Ž  des permutations qui
laissent fixe le point $e$. On v\Žrifie, sans grand
d\Žtour, que la structure $T/R$ est isomorphe \ˆ la
structure $G/S$ elle-m\me isomorphe au D-hypergroupe
$G/H$.\cqfd

\

On trouvera davantage de pr\Žcisions au sujet de ces
structures dans l'\sc Appendice\rm, ci-dessous.

\

Pour le moment, on notera, cependant, encore ceci : pour

$I =
\{0,1\} \ , \ |A_0| = n \geqs 3 \ , \ |A_1| \geqs 3 \ , \ n
\neq p$, 

\noi on obtient une structure \ˆ 
$n+p$ \Žl\Žments qui est un hypergroupe, mais {\bf pas un}
D-hypergroupe : voir le {\bf A5} de l'\sc
Appendice\rm.

\

\head{Scholie} 

\

On va explorer les liens qui peuvent exister entre deux pr\Žsentations ayant une m\me trame
$(T,C)$.

\

On se donne deux relations d'\Žquivalence $R$ et $S$ sur
l'ensemble $T$. On d\Žsigne, respectivement, par $\bar x$ et
$\hat x$ la classe de $x$ modulo $R$ et la classe de $x$ modulo $S$. On a
ainsi deux pr\Žsentations, $(T,C,R)$ et $(T,C,S)$, de la structure quotient
$T/R$ et de la structure quotient $T/S$, respectivement,  qui sont munies des
op\Žrations multivalentes correspondantes : 
$$\bar z \in \bar x .\bar y \iff (\exists (u,v) \in (\bar x \times \bar y) \cap
C) ( uv \in \bar z),$$
$$\hat z \in \hat x .\hat y \iff (\exists (u,v) \in (\hat x \times \hat y) \cap
C)( uv \in \hat z).$$

\

Bien entendu, lorsque
$R \inc S$, on a $\bar x \subset \hat x$ et, dans ce cas, 
$\bar x \mapsto \hat x$ est une application surjective
$f : T/R
\to T/S$, $f(\bar x) = \hat x$. De plus, il suffit de se
reporter aux d\Žfinitions pour constater que l'on aura
alors, toujours, l'inclusion

$f(\bar x.\bar y) \inc \hat x.\hat y$.

\noi Pour avoir l'\Žgalit\Ž

$f(\bar x.\bar y) = \hat x.\hat y$, 

\noi il faut et il suffit
que l'on ait
 
$(**)((u,v) \in (\hat x \times \hat y) \cap C)
\implies$

$ (\forall (p,q) \in \hat x
\times \hat y)(\exists (r,s) \in (\bar p
\times \bar q) \cap C)(\widehat{rs} =
\widehat{uv})$.

\

Pour examiner plus attentivement cette situation, commen\c
cons par lui donner un nom commode :  convenons de dire que
l'\Žquivalence $S$ est
{\bf invariante modulo}
$R$ lorsque l'on a $R \inc S$ et que la condition pr\Žc\Ždente, ($**$),
est satisfaite pour tous $x,y$ dans $T$.

\

\bf Supposons donc que $\bold S$ soit invariante modulo
$\bold R$. \rm Il va en d\Žcouler {\bf trois} cons\Žquences importantes.

\su{1} Pour tous $x,y,z$ dans $T$, on aura
$$f^{-1}f(\bar x.\bar y) = f^{-1}(f(\bar x).f(\bar y)) =
\bar x.f^{-1} f(\bar y) =f^{-1}f(\bar x).\bar y.$$

\su{2} Si $R$ est ad\Žquate, il en sera de m\me de $S$.
Autrement dit, si $T/R$ est un hypergroupe, $T/S$ sera un
reflet de $T/R$ car $f$ est alors un r\Žflecteur d'apr\s
{\bf 1}.

\su{3} En particulier, si $T/R$ est un groupe, $T/S$ sera une
image homomorphe de $T/S$.

\su{D\Žmonstration}

\su{1}  Du fait de l'\Žgalit\Ž 
$$f(\bar x.\bar y) = \hat x.\hat y = f(\bar x).f(\bar y),$$
on a aussi
$$f^{-1}f(\bar x.\bar y) = f^{-1}(f(\bar x).f(\bar y)).$$
Par d\Žfinition, on a
$$f(\bar x) = \hat x \ , \ f(\bar y) = \hat y.$$ Il nous
faut \Žtablir l'identit\Ž 
$$f^{-1}(f(\bar x).f(\bar y)) = \bar x.f^{-1} f(\bar y)$$
qui s'\Žcrit donc
$$f^{-1}(\hat x.\hat y) = \bar x.f^{-1}(\hat y).$$ 
On proc\de par implications et \Žquivalences successives. On a
$$\bar z \in f^{-1}(\hat x.\hat y) \Longleftrightarrow \hat
z \in \hat x.\hat y \Longleftrightarrow \exists \ (u,v) \in
(\hat x \times \hat y) \cap C \ \ \text{tel que} \ \  uv \in
\hat z,$$
$$\bar z \in \bar x.f^{-1}(\hat y) \Longleftrightarrow
\exists \ t \in \hat y \ \ \text{et} \ \ \exists \ (u,v) \in
(\bar x
\times
\bar t) \cap C \ \ \text{tel que} \ \ uv \in \bar z.$$ Or,
$$t \in \hat y \ \ \text{et} \ \ (u,v) \in (\hat x
\times
\hat t) \cap C 
\Longrightarrow (u,v) \in (\hat x
\times
\hat y) \cap C  ,$$ donc
$$\bar z \in \bar x.f^{-1}(\hat y) \Longrightarrow  \bar z
\in f^{-1}(\hat x.\hat y).$$ 
Pour la r\Žciproque, on
utilise la condition ($**$) sous la forme suivante :
$$t \in \hat y \ \ \text{et} \ \ (u,v) \in (\hat x \times
\hat y) \cap C \Longrightarrow \exists \ (r,s) \in (\bar x
\times \bar t) \cap C \ \ \text{tel que} \ \  \widehat{rs} =
\widehat{uv},$$ donc
$$ \bar z
\in f^{-1}(\hat x.\hat y) \Longrightarrow \bar z \in \bar
x.f^{-1}(\hat y).$$ D'o\  l'\Žgalit\Ž annonc\Že
$$f^{-1}(f(\bar x).f(\bar y)) = \bar x.f^{-1} f(\bar y).$$
On \Žtalit de m\me, l'\Žgalit\Ž 
$$f^{-1}(f(\bar x).f(\bar y)) = f^{-1}f(\bar x).\bar y.$$
Cela d\Žmontre le point {\bf 1}.

\su{2} Si $T/R$ est un hypergroupe, il en va de m\me de
$T/S$ du fait que $f$ est surjective et v\Žrifie
l'identit\Ž  
$$f(\bar x.\bar y) = \hat x.\hat y = f(\bar x).f(\bar y).$$
\se, l'op\Žration multivalente de $T/S$ est
alors {\bf associative} car celle de $T/R$ l'est :
$$(\hat x.\hat y).\hat z = f(\bar x.\bar y).f(\bar z) =
f(\bar x.\bar y.\bar z) = f(\bar x).f(\bar y.\bar z) = \hat
x.(\hat y.\hat z).$$ 

De m\me, elle poss\de \Žgalement
la propri\Žt\Ž de {\bf reproductibilit\Ž}  comme pour
$T/R$.

\

Cela d\Žmontre le point {\bf 2} dont d\Žcoule imm\Ždiatement
le point {\bf 3}.\qed

\

On aura remarqu\Ž, sans doute, la similitude avec la d\Žmonstration faite au {\bf A1} de l'\sc Appendice \rm pour
le cas particulier des D-hypergroupes.

\

On en vient, \ˆ pr\Žsent, \ˆ la g\Žn\Žralisation
compl\te suivante du th\Žor\me qui caract\Žrise les
reflets des D-hypergroupes.

\su{Th\Žor\me} \sl Soit $(T,C,R)$ une pr\Žsentation d'un
hypergroupe $H$ quelconque. Tout reflet $L$ de $H$ poss\de
une pr\Žsentation $(T,C,S)$ o\ $S$ est une
\Žquivalence invariante modulo $R$\rm.

\su{D\Žmonstration} Soient $p : T \to h = T/R$ la
projection de $T$ sur l'ensemble quotient, $h : H
\to L$ un r\Žflecteur et $S$ l'\Žquivalence sur $T$ d\Žfinie par l'application compos\Že $g = hp$, autrement dit,
$$x S y \iff h(p(x)) = h(p(y)).$$
L'hypergroupe $L$ est ainsi identifi\Ž \ˆ  l'hypergroupe
$T/S$ pr\Žsent\Ž par $(T,C,S)$. On adopte les notations
pr\Žc\Ždentes, en posant
$$\bar x = p(x) \ , \ \hat x = h(p(x)).$$

On va montrer que $S$ v\Žrifie la condition ($**$). On se
donne

$(u,v) \in (\hat x \times \hat y) \cap C$ et $(p,q) \in \hat
x \times \hat y$. 

On a ainsi 

$\hat p = \hat x$, $\hat q =
\hat y$ et
$\widehat{uv} \in \hat x .\hat y = h(\bar p.\bar q)$.

Puisque $h$ est un r\Žflecteur, il existe $z$ tel que
$\bar z
\in
\bar p.\bar q$ et $\hat z =  \widehat{uv}$. Cela 

veut dire
qu'il existe $(r,s) \in (\bar p \times \bar q) \cap C$ tel
que $\bar z = \overline{rs}$ et, par suite,

$\widehat{rs} = \hat z  =  \widehat{uv}$.\qed

\su{Corollaire} \sl Soit $(T,C,R)$ une pr\Žsentation d'un
hypergroupe $H$ non trivial. Cet hypergroupe $H$ est
{\bf simple} si et seulement si les deux seules
\Žquivalences invariantes modulo
$R$ sont l'\Žquivalence totale $T \times T$ et $R$
elle-m\me\rm.

\su{Une remarque} Insistons : 

\

{\bf tout hypergroupe poss\de une pr\Žsentation} \ˆ 
l'aide d'une trame convenable.

\

La d\Žmonstration donn\Že plus haut \Žtablit, en fait, le
r\Žsultat encore plus g\Žn\Žral que voici. Soit $H$ un
ensemble muni d'une op\Žration multivalente quelconque :
autrement dit, certains des produits $x.y$ peuvent \tre
{\bf vides} et l'op\Žration n'a nul besoin d'\tre
associative, ni reproductive. La construction donn\Že plus
haut fournit une trame $(T,C)$ et une relation d'\Žquivalence $R$ sur $T$ telles que $(T,C,R)$ soit une pr\Žsentation de $H$, autrement dit, $H$ est isomorphe \ˆ  la
structure quotient $T/R$.

\

On peut concevoir cette pr\Žsentation comme un {\bf d\Žploiement} de l'op\Žration multivalente sur $H$ en une
op\Žration univalente sur $T$, mais partiellement d\Žfinie.
Ce que l'on gagne en simplicit\Ž, en passant du multivalent
\ˆ
l'univalent, on le perd, un tout petit peu, en passant
d'une op\Žration partout d\Žfinie \ˆ  une op\Žration
uniquement d\Žfinie pour certains couples : le mal n'est
pas si grand car cela sous-tend la notion d'une op\Žration
univalente {\it g\Žn\Žrale} pour laquelle certains
produits
$x.y$ au lieu d'\tre des singletons, seraient simplement
{\bf vides}. Dans ce d\Žploiement, chaque \Žl\Žment $x$
de
$H$ est, pour ainsi dire, {\bf \Žclat\Ž} en une partie
de
$T$.

\

On peut aussi concevoir cette pr\Žsentation comme un {\bf
rev\tement} de la structure $H$ par la structure
$(T,C)$, un peu comme pour les vari\Žt\Žs topologiques.

\

Cette fa\c con de faire est loin d'\tre rare. Elle est
illustr\Že par de tr\s nombreux exemples. Citons-en
uniquement deux : le rev\tement d'un groupe de Lie connexe
par un groupe de Lie simplement connexe, et la repr\Žsentation d'un ensemble analytique de la droite r\Želle 
par la projection d'un ferm\Ž du plan. Les chemins qui ne
se distinguent pas tr\s bien dans le groupe, se voient
bien plus clairement une fois relev\Žs dans le rev\tement.
De m\me, le sch\Žma (complexe) de Souslin qui d\Žfinit
l'ensemble analytique est d\Žploy\Ž, \ˆ  l'aide du ferm\Ž
du plan, comme une nappe froiss\Že et entortill\Že que
l'on \Žtend sur une corde \ˆ  linge.

\

Bien entendu, le crit\re de simplicit\Ž  qu'on vient de
donner pour les hypergoupes s'applique \ˆ  toutes les
structures multivalentes ou univalentes. Il s'applique, en
particulier, aux groupes : un groupe est simple si et
seulement si l'une de ses pr\Žsentations $(T,C,R)$, v\Žrifie la condition suivante :  les deux seules
\Žquivalences invariantes modulo
$R$ sont l'\Žquivalence totale $T \times T$  et $R$
elle-m\me\rm. De plus, s'il en est ainsi, alors {\bf toutes} les pr\Žsentations du groupe v\Žrifie cette m\me
condition.

\

La question se pose alors : est-il plus simple d'\Žtablir
qu'un groupe est simple, directement, en s'y tenant au plus
pr\s, ou bien peut-on simplifier la preuve, parfois, en
prenant de la hauteur au-dessus du groupe, \ˆ  l'aide d'une
trame convenable ?

\

\head{Constructions \guil \emph{\ˆ la Utumi}\guir}

\

\

On se donne un hypergroupe $H$ dont on d\Žsigne l'op\Žration par le signe $+$ (m\me si elle n'est pas
nŽcessairement commutative). On se donne une partition
$\cal P$ de l'ensemble $H$ et, pour chaque $x \in H$, on d\Žsigne par
$\bar x$ la partie de la partition $\cal P$ \ˆ laquelle $x$
appartient (sa \guil classe\guir). On suppose que
$H$ poss\de un \Žl\Žment $0$ ayant les propri\Žt\Žs suivantes : 
$$\bar 0 = \{0\} \ , \ x + 0 = x \ \ \text{et} \ \  x \in 0 +
x 
\inc
\bar x
\pourtout
\ x
\in H.$$ 
On notera que l'on a alors $0 + \bar x = \bar x$.

\

\noi On d\Žfinit sur $H$ une nouvelle op\Žration
multivalente par la r\gle :
$$x.y = x + \bar y.$$ 
Cette op\Žration est toujours
reproductive puisque
$$x.H = x + H = H \et \ H.x = H + \bar x = H.$$
Elle est associative si et seulement si, pour tous $x$ et
$y$, on a
$$\bar x + \bar y = \overline{x + \bar y}.$$ 
Cela se voit
simplement, sans grand d\Žtour. On dira alors que cet hypergroupe $(H,.)$ est
construit
\guil \ˆ la Utumi\guir \  \ˆ partir de l'hypergroupe $(H,+)$
et de la partition $\cal P$.

\su{L'exemple de Utumi} 

\ 

Voici le premier exemple de cette construction, donn\Ž en 1949 par Utumi [7]. On
prend pour $H$ le groupe cyclique d'ordre huit :

$H = \Z/8 \Z =  \{0,1,2,3,4,5,6,7\}$

et la partition form\Že des
trois classes :

$I = \{0\}$, $A = \{1,4,7\}$ et
$B =
\{2,3,5,6\}$. 

\noi Les conditions pr\Žc\Ždentes, \Ždict\Žes pour la construction, \Žtant
satisfaite, $(H,.)$ est un hypergroupe. Utumi  a
montr\Ž
que ce n'\Žtait pas un D-hypergroupe bien qu'il poss\de
la propri\Žt\Ž  suivante : pour chaque \Žl\Žment fix\Ž
\ $x
\in H$, la collection
$\{x.y : y \in H \}$ est une partition de $H$ et les
cardinaux
$|x.y|$ sont tous
\Žgaux, pour
$y
\in H$.

En 1940, Eaton  [1] a donn\Ž le nom de {\it
cogroupes} (\ˆ
droite) \ˆ cette classe particuli\re d'hypergroupes et
l'hypergroupe $H$ de Utumi en fait partie. Aussi lui a-t-on
donn\Ž
le nom de
{\it cogroupe de Utumi}. Pour plus de d\Žtail, on
se reportera
\ˆ
[3]. 

On va  \Žtablir, ici, le r\Žsultat
nouveau que voici.

\su{Th\Žor\me} \sl Le cogroupe de Utumi est un hypergroupe
simple\rm.

\ 

\se, soit $f : H \to K$ un r\Žflecteur et
pour chaque partie $X \inc H$, posons $\hat X =
f^{-1}f(X)$ et, en particulier, $\hat x = f^{-1}f(x)$ pour
chaque $x \in H$. Les  propri\Žt\Žs du  r\Žflecteur se
traduisent, ici, de la mani\re suivante :
$$\widehat{x.y} = \hat x.\hat y = x.\hat y = \hat x . y.$$
De ce que $x.0 = x$, on tire :
$$\hat x = \hat x.\hat 0 = x.\hat 0 = \hat x0$$
et, en particulier,
$$\hat 0 = \hat 0.\hat 0 = 0.\hat 0 = \hat 0.0.$$
Deux cas seulement peuvent se pr\Žsenter.

\

Premier cas : $\hat 0 = 0$. Alors, pour chaque $x$, on a
$\hat x = x.0 = x$ et $f$ est un isomorphisme.

\  

Second cas : $\hat 0$ contient un \Žl\Žment $a \neq 0$.
Alors $\hat 0$ contient aussi le  produit $0.a =  \bar a$
donc,
\Žgalement, le produit $\bar a.a = \bar a + \bar a$ ainsi
que toutes les sommes $\bar a + \bar a + \dots + \bar a$. Une
simple v\Žrification montre que $A + A + A = H$ et $B + B =
H$. Donc, quel que soit $a \neq 0$, on voit que l'on a $\hat
0 = H$ et l'hypergoupe $K$ est donc trivial, r\Žduit \ˆ un
singleton.\cqfd

\

Plus g\Žn\Žralement, et en brodant sur ce m\me canevas,
on  obtient donc ceci.

\

\sl Un hypergroupe $(H,.)$ construit
\guil \ˆ la Utumi\guir,  \ˆ partir d'un groupe $(H,+)$ et
d'une partition $\cal P$, est {\bf simple} \sl d\s que, pour  chacune des parties
$A
\in
\cal P$, diff\Žrente du singleton $\{0\}$, l'une au moins des
sommes
$A + A + \dots + A$ est \Žgale \ˆ $H$\rm.

\

\head{Appendice}

\

Lorsque l'on traite d'hypergroupes, on adopte souvent une
convention tacite : chaque fois que le produit $x.y$ est un
singleton, on \Žcrit, indiff\Žremment, $x.y =
\{z\}$ ou $x.y = z$, pourvu que le risque de confusion soit
minime.

\

\head{Reflets d'un D-hypergroupe}

\

\su{A1. Lemme} \sl Soient $H$ et $K$ deux sous-groupes d'un
groupe
$G$. On suppose que $K$ est invariant modulo $H$. Le
D-hypergroupe $G/K$ est alors un reflet du D-hypergoupe
$G/H$\rm.

\

\su{D\Žmonstration} Pour $x \in G$, posons $\bar x = xH$.
Par hypoth\se,
$H$ est un sous-groupe de
$K$. En associant \ˆ  chaque classe $ \bar x = xH \in G/H$
la classe
$\hat x = xK \in G/K$, on obtient donc une application
surjective canonique
$f : G/H
\to G/K$. On va v\Žrifier que $f$
est un r\Žflecteur.

On a
$f(\bar x) = \hat x =  xK$, de sorte que
$$f(\bar x . \bar y) = \{ \hat z : z \in xHyH \},$$
$$f(\bar x).f(\bar y) = \hat x.\hat y = \{ \hat z : z \in
xKyK
\},$$ et ces deux ensembles sont
\Žgaux. 

En effet, soit  $z \in xKyK$ alors $z \in xHyK$,
puisque $KyK = HyK$; autrement dit, l'intersection
$xK \cap (xHy)$ n'est pas vide; soit $t$ un \Žl\Žment de
cette intersection; alors $t \in xHyH$ et $zK = tK$. On a
ainsi
$$f(\bar x . \bar y) = f(\bar x).f(\bar y).$$ 
De m\me, on a 

$f^{-1}f(\bar x) = f^{-1}(\hat x) = \{ \bar z : z \in xK
\}$, 

de sorte que

$$f^{-1}f(\bar x.\bar y) = f^{-1}(f(\bar x).f(\bar y)) = \{
\bar z : z \in xKyK \},$$
$$f^{-1}f(\bar x).\bar y = \{ \bar z : z \in xKyH \},$$
$$\bar x.f^{-1}f(\bar y) = \{ \bar z : z \in xHyK \}.$$ 
Ces
quatre ensembles sont \Žgaux et l'application
$f$ est donc bien un r\Žflecteur.\qed

\su{A2. Deux remarques utiles} Pour \Žtablir le lemme
suivant, on aura besoin des deux r\Žsultats que voici.

\su{R1} Tout r\Žflecteur $f$ v\Žrifie l'identit\Ž 
$$f^{-1}f(x).f^{-1}f(y) = f^{-1}(f(x).f(y)).$$
\se, soit $z \in f^{-1}f(x).f^{-1}f(y)$. Il
existe alors deux \Žl\Žments

$u \in f^{-1}f(x)$ et $v \in f^{-1}f(y)$ tels que $z \in u.v$. 

On a donc $f(z) \in f(u.v) = f(u).f(v) = f(x).f(y)$.

Autrement dit, on a $f^{-1}f(x).f^{-1}f(y) \inc
f^{-1}(f(x).f(y))$. 

On conclut en remarquant que
$$f^{-1}(f(x).f(y)) = x.f^{-1}f(y) \inc
f^{-1}f(x).f^{-1}f(y). \ \cqfd$$

\

\su{R2} On se donne des hypergoupes,  $H,K,L,$ ainsi que des
applications, $f : H
\to K \ , \ g : K \to L$, et leur compos\Že $h = gf$.  Si $f$
et $g$ sont des morphismes, il est clair que $h$ est \Žgalement un morphisme. R\Žciproquement, si $f$ et $h$ sont
des morphismes, et pourvu que $f$ soit
{\bf surjective},  l'application $g$ est aussi un morphisme.

\

\se, pour $x,y$ dans $K$, il existe $a,b$
dans $H$ tels que $x = f(a)$ et $y = f(b)$. Ainsi $f(a.b) =
f(a).f(b) = x.y$ donc 

$g(x.y) = h(a.b) = h(a).h(b) =
gf(a).gf(b) = g(x).g(y)$. 

\noi Autrement dit, $g$ est  bien un
morphisme.\cqfd

\su{A3. Lemme} \sl Soient $G$ un groupe, $H$ l'un quelconque
de ses sous-groupes et $L$ un reflet du D-hyprgroupe $G/H$.
Il existe alors un sous-groupe $K$ de $G$, invariant modulo
$H$, tel que $L$ soit isomorphe au D-hypergroupe $G/K$\rm.

\su{D\Žmonstration} Soit
$h : G/H
\to L$  un r\Žflecteur. D\Žsignons par
$e$ l'\Žl\Žment unit\Ž
du goupe G puis, pour chaque $x \in G$, posons  $\bar x =
xH$. Soient $1 = h(\bar e)$ et
$K = \{x : h(\bar x) = 1 \}$. On proc\de par \Žtapes. On
commence par montrer que
$K$ est  un sous-groupe de $G$ puis on montre qu'il est
invariant modulo
$H$ et on \Žtablit, enfin, que $L$ est isomorphe \ˆ 
l'hypergroupe
$G/K$.

\su{1} Pour $x \in H$, on a $\bar x = \bar e$, de sorte que
$h(\bar x) = h(\bar e) = 1$. On a donc
$$H \inc K.$$
 
On commence par utiliser le fait que $h$ est un morphisme.

\su{2} Pour
$x
\in G$ et
$y
\in G$, on a $xy \in xHyH$ donc $\overline{xy} \in (xH).(yH)
= 
\bar x .
\bar y$, de sorte que
$$h(\overline{xy}) \in h(\bar x . \bar y) = h(\bar x) .
h(\bar y).$$

\su{3} Pour $x \in G$, on a aussi $\bar x . \bar e = \bar x$
d'o\
l'on tire $h(\bar x).1 = h(\bar x)$. De sorte que, pour $y
\in K$, on a
$$h(\overline{xy}) \in h(\bar x).1 = h(\bar x) \ \
\text{qui est un singleton, donc} \ \ h(\overline{xy}) =
h(\bar x).$$

\su{4} En particulier, pour $y \in K$, il vient
$$1 = h(\overline{y^{-1}y}) =  h(\overline{y^{-1}}) \donc
\ y^{-1} \in K.$$

\su{5} Plus g\Žn\Žralement, si $x \in K$ et $y \in K$,
alors
$$h(\overline{x^{-1}y)} = h(\overline{x^{-1}}) = 1 \donc
\ x^{-1}y
\in K,$$ 
de sorte que 
$$K^{-1}K \inc K.$$
{\bf Cela prouve bien que $K$ est un sous-groupe de
$G$}.

\su{6} En proc\Ždant par \Žquivalences successives, on a
$$\bar z \in h^{-1}h(\bar e) \iff \bar z \in
h^{-1}(1) 
\iff h(\bar z) = 1 \iff
z \in K.$$
Autrement dit
$$\bar z \in h^{-1}h(\bar e) \iff z \in K.$$
On utilise, \ˆ  pr\Žsent, les identit\Žs sp\Žcifiques aux
r\Žflecteurs.

\su{7} On a $\h h(\bar x) = \h h(\bar x . \bar e) = \bar x
.\h h(\bar e)$. D'o\ 

$\bar z \in \h h(\bar x) \iff( \exists \bar u \in \h h(\bar e))(\bar z \in \bar x
. \bar u) \iff$

$(\exists u \in K)(z \in xHuH)
\iff z \in xK$. 

Ainsi, $h(\bar z) = h(\bar x) \iff z \in
xK$. Il en r\Žsulte que l'on a 
$$h(\bar x) = h(\bar y) \iff xK = yK.$$ Il
existe ainsi une application $g : G/K \to L$ qui, \ˆ 
chacune des classes $xK \in G/K$, associe $g(xK) = h(\bar x) \in L$
et $g$ est une application {\bf bijective}.

\su{8} On a l'identit\Ž  $\bar e. \h h(\bar x) = \h h(\bar e). \h h(\bar x)$, d'apr\s {\bf R1}. 

Or, $\bar z \in \bar
e.
\h h(\bar x)
\iff z \in HxK$. 

De m\me, $\bar z \in \h h(\bar
e). \h h(\bar x) \iff z \in KxK$. 

Donc, pour tout $x \in
G$, on a 
$$HxK = KxK.$$ 
{\bf Cela ach\ve de prouver que $K$ est invariant modulo
$H$}.

\su{9} Soit $f : G/H \to G/K$ l'application canonique qui
fait correspondre, \ˆ  la classe $xH \in G/H$, la classe
$xK \in G/K$. C'est un r\Žflecteur, d'apr\s {\bf A1}. On
a
$g(f(\bar x)) = g(xK) = h(\bar x)$. Autrement dit, $h$ est la
compos\Že $gf$. Or, en particulier, $f$ est un morphisme
surjectif et 
$h$ un morphisme. Donc $g$ est un morphisme, d'apr\s {\bf
R2}. Puisque
$g$ est bijectif (voir le {\bf 7} ci-dessus) c'est donc un
isomorphisme, comme annonc\Ž.\qed

\su{A4. Conditions d'ad\Žquation} Soient $(T,C)$ une trame et $R$ une relation d'\Žquivalence sur $T$. Pour chaque $x \in T$, on d\Žsignera par $\bar x$
sa classe modulo $R$, de sorte que toutes les relations
suivantes sont synonymes :
$$x R y \ , \ y R x \ , \ y \in \bar x \ , \ x \in \bar y \ ,
\ \bar x = \bar y.$$ 
Afin de simplifier les \Žcritures,
pour $r \in T$ et $s \in T$, on conviendra de ceci : lorsque
l'on \Žcrit $rs$, on suppose (implicitement) que le couple
$(r,s)$ est composable, autrement dit, que l'on a $(r,s) \in
C$.

\

Voici deux \Žnonc\Žs  relatifs \ˆ  la
structure quotient $T/R$, (ils sont tous deux du premier ordre) :
\[(\forall x \ \forall y) (\exists \ r \in \bar x \
\exists
\ s
\in T)
(rs \in \bar y) \et 
(\exists \ r \in T \ \exists  \ s \in \bar y)
( rs \in \bar x)\}.\tag{1}\]
\[(\forall x)(\forall y)(\forall z)(\forall u)(\exists r
\in \bar x \ , \ s \in \bar y \ , \ t \in \bar z)((rs)t
\in \bar u\ \iff \tag{2}\]
\[(\exists r
\in \bar x \ , \ s \in \bar y \ , \ t \in \bar z)(r(st)
\in \bar u).\]

\

On voit, sans d\Žtour, que la condition (1) est n\Žcessaire et suffisante pour assurer la {\bf reproductibilit\Ž}. De m\me, la condition (2) est n\Žcessaire et
suffisante pour assurer {\bf l' associativit\Ž}.

\

Autement dit, pour que l'\Žquivalence $R$ soit {\bf ad\Žquate}, il faut et il suffit qu'elle satisfasse les deux \Žnonc\Ž s (1) et (2).

\

\head{Sur la nouvelle famille d'exemples}

\

On se donne un ensemble $K$, r\Žunion d'une famille
$(A_i)_{i \in I}$ d'ensembles non vides, deux \ˆ deux
disjonts. On suppose que
$0
\in I$ et on distingue un point particulier
$e
\in A_0$. 

\

Sur cet ensemble
$K$, on a d\Žfini une op\Žration multivalente dont voici la
table : quels que soient $a_i \in A_i$, pour chaque indice $i$,

$x.e = x, \pourtout  \ x \in K$,

$a_i.y = A_i \stm \{a_i\}, \pourtout  \ y
\in A_0 \stm \{e\}$,

$a_i.a_j = K \stm A_i, \pourtout \ \text{indice} \
j \neq 0$.

\

Il est clair que la structure ainsi d\Žfinie sur
$K$ d\Žpend, uniquement, \ˆ
isomorphisme pr\s, de la famille des cardinaux,  finis ou infinis, $n =
|A_0|$ et, pour $i \neq 0$, $p_i = |A_i|$ : on dira qu'elle
est de type $S(n,(p_i))$. Cela g\Žn\Žralise, quelque peu,
les structures de type
$S(n,p)$, correspondant au cas o\ l'on a
$|I| \leqs 2$, introduite dans l'article [2]. 

\su{A5}  Lorsque l'on a $n \geqs 3$, $p_i \geqs 3$, pour tout indice
$i$, et  $p_k \neq n$, pour un indice $k$, la structure est un
hypergroupe qui n'est pas un D-hypergroupe.

\

\se, pour voir que la structure est un
hypergroupe, il suffit de s'assurer que les conditions
d'ad\Žquation sont satisfaites.

\

Pour la condition de reproductibilit\Ž , la v\Žrification est
sans d\Žtour.

\

Quant \ˆ  la condition d'associativit\Ž, on reprend la
trame $(T,C)$ qui sert \ˆ construire la structure. Soient
alors
$a_i,a_j,a_k,$ donn\Žs dans
$K$ ainsi que $r,s,t,$ dans $T$ tels que $r(e) = a_i$,
$s(e) \in a_j$, $t(e) \in a_k$. Si $(r,s) \in C$ et $(rs,t)
\in C$, on voit, simplement, que l'on a $(s,t) \in C$ et
$(r,st)
\in C$, de sorte que
$((rs)t)(e) = (r(st))(e)$. 

\

R\Žciproquement, soit $(s,t) \in C$ et $(r,st) \in C$.
Ainsi,
$s(t(e)) = s(a_k)$ est d\Žfini, de m\me que $r(s(a_k)) =
a_l$. On montre qu'il existe alors un \Žl\Žment $q \in T$
d\Žfini au point $a_j$ et tel que 
$$q(e) = r(e) = a_i \ , \ (q,s) \in C \ , \ q(s(a_k)) = 
r(s(a_k)) = a_l \ ,
\ (qs,t)
\in C,$$ 
ce qui \Žtablira la r\Žciproque. L'application
$r$ est d\Žj\ˆ d\Žfinie au point $e$ et au point $s(a_k)$. Si $a_j$ est
l'un des points
$e$ ou
$s(a_k)$, on prend $q = r$. Sinon, on construit $q :
\{e,s(a_k),a_j\}
\to K$ en posant 
$$q(e) = r(e) = a_i \ , \ q(s(a_k)) = r(s(a_k)) = a_l \ , \
q(a_j) = x,$$  avec un \Žl\Žment $x$ soumis aux trois conditions
suivantes :

(1) $x
\notin \{a_i,a_l\}$.

(2) $x
\in A_i$ si $e \in A_j$. 

(3) $x \in A_l$ si $s(a_k) \in A_j$.

\noi Cela est possible car chacune des parties 
$A_m$ de la famille poss\de au moins trois points. Enfin,
ces conditions impliquent que l'on a bien $q \in T$.

\

On pourrait \Žgalement  v\Žrifier \ˆ l'aide de la table
de multiplication, avec un peu de patience, que l'op\Žration est associative et reproductive.

\

Enfin, pour $y \in A_0 \stm \{e\}$ et $a \in A_k$, les
deux parties $e.y = A_0 \stm \{e\}$ et $a.y = A_k
\stm \{a\}$ ne sont pas \Žquipotentes alors qu'elles
le seraient dans un D-hypergroupe.\cqfd

\

On peut ainsi ranger les diff\Žrentes structures
$S(n,(p_i))$ dans quatre classes distinctes.

\

Premi\re classe : $p_i = n$ pour tout indice
$i$.  Cela comprend, bien entendu, tous les cas o\ 
$I = \{0\}$. Toutes ces structures sont des D-hypergroupes,
comme on l'a dit plus haut.

\

Deuxi\me classe : $n \geqs 3$, $p_i \geqs 3$ pour tout indice
$i$ et  $p_k \neq n$ pour un indice $k$. Ces structures sont,
toutes, des hypergroupes qui ne sont pas des hypergroupes de
classes.

\

Troisi\me classe : $n \geqs 2$  et $p_k = 1$ pour un indice
$k$. L'un des produits est {\bf vide}  : en particulier,
$A_0 = \{e,y,\dots\}$,  $A_k =
\{a\}$ et 
$a.y =
\vide$.

\

Quatri\me classe : tous les autres cas. Ces structures ne
sont pas des hypergroupes car l'associativit\Ž est alors en d\Žfaut, comme on peut le v\Žrifier.

\

\su{V\Žrifications} Il n'y a que les trois cas suivants :

\su{1} $n = 1$ et $p_k \geqs 2$ pour un indice $k$.  

\su{2} $n = 2$ et $p_k \geqs 3$ pour un indice $k$.

\su{3} $n \geqs 3$ et $p_k = 2$ pour un indice $k$.

\su{1} On a $A_0 = \{e\}$, $A_k = \{a,b,\dots\}$ donc $a.(a.a)
\neq (a.a).a$ car
$$a.(a.a) = a. (K \setminus A_k) \subset \{a\} \cup (K
\setminus A_k) \not\ni b,$$
$$(a.a).a = (K \setminus A_k).a = \bigcup_{i \neq k} (K
\setminus A_i) \ni b.$$

\su{2} On a $A_0 = \{e,y\}$, $A_k = \{a,b,c,\dots\}$, donc
$a.(y.y) \neq (a.y).y$ car
$$a.(y.y) = a.e = a,$$
$$(a.y).y = (A_k \setminus \{a\}).y \supset (b.y) \cup (c.y)
= A_k.$$

\su{3} On a $A_0 = \{e,y,z\}$, $A_k = \{a,b\}$ donc $a.(y.y)
\neq (a.y).y$ car
$$a.(y.y) = a.\{e,z\} = A_k,$$
$$(a.y).y = b.y = a.\qed$$

\

\head{Pour la petite histoire}

\

Les premi\res versions de ce texte remontent au si\cle dernier (!) Tr\s souvent remani\Ž, il a \Žt\Ž expos\Ž, en partie, 
au Groupe de travail de th\Žorie des nombres du D\Žpartement de math\Žmatiques de l'Universit\Ž Blaise Pascal, de Clermont-Ferrand, sur invitation, le lundi 18 juin 2012, sous le titre :

\guil Est-il encore possible d'avancer \ˆ la mani\re d'Emile Mathieu ?\guir

\noi avec le r\Žsum\Ž suivant : 

\guil Courte promenade math\Žmatique \ˆ travers champs, suivie d'un bref plaidoyer en faveur de la r\Žvision de la liste des sporadiques.\guir

\

\head{En guise de conclusion (provisoire)}

\

On sait faire agir les hypergroupes sur des ensembles, comme on le fait pour les groupes. Plus pr\Žcis\Žment, on sait les faire agir sur l'ensemble des racines d'un polyn\™me. Et c'est encore une autre histoire. On peut d\s lors se poser la question suivante. Est-il raisonnable, ou non, de pr\Ždire l'existence d'hypergroupes de Mathieu simples ?

\

\guil \`A peine clos, le chantier de la classification des
{\bf groupes} simples finis devra-t-il, \Žventuellement,  rouvrir ses portes
afin d'entreprendre la classification des {\bf hypergroupes} simples finis?\guir

\head{Bibliographie}

\

\noi [1] Eaton, J. E., Theory of cogroups, \it Duke
Math. J.\rm, {\bf}  6 (1940) 101-107.

\

\noi [2] Haddad, Labib  et Sureau, Yves,  Les cogroupes
et les
$D$-hypergroupes, \it J.
Algebra\rm, {\bf 118} (1988) 468-476.

\

\noi [3] Haddad, Labib et Sureau, Yves,  Les cogroupes
et la construction de Utumi, \it Pacific J. Math.\rm, {\bf 145} (1990) 17-58.

\

\noi [4] Krasner, M., Sur la th\Žorie de la ramification des id\Žaux de
corps nongaloisiens de nombres alg\Žbriques,  \it Acad.
Belgique, Cl. Sci. M\Žm. Coll.  4$^\circ$\rm, {\bf 11}, 
Th\se,  (1937), 110 pages.

\

\noi [5] Marty, F., Sur une g\Žn\Žralisation de la notion de groupe, \it 8th. Scand. Math. Congr.\rm,  Stockholm (1934) 45-49.

\

\noi [6] Marty, F., Sur les groupes et hypergroupes attach\Žs \ˆ une
fraction rationnelle, \it Ann. Sci. \'Ec. Norm.
Sup\Žr.  Ser. III\rm. {\bf 53} (1936)  83-123.

\

\noi [7] Utumi, Yuzo, On hypergroups of group right
cosets, \it Osaka Math. J.\rm, {\bf 1} (1949) 73-80.

\

\

\

\
 
\enddocument